\documentclass[11pt]{article}
\title{A note on special duality triads and their operator valued counterparts}
\author{Ewa Borak \\
\\Institute of Computer Science, Bia{\l}ystok University\\
PL-15-887 Bia{\l}ystok, ul.Sosnowa 64, POLAND\\
e-mail: ewag@ii.uwb.edu.pl, ewaborak@wp.pl}

\usepackage{amsmath,amsthm}

\chardef\bslash=`\\ 
\newtheorem{ex}{Example}[section]
\newtheorem{defn}{Definition}[section]
\newtheorem{obs}{Observation}[section]
\newtheorem{cor}{Corollary}[section]

\newtheorem{prop}{Proposition}[section]
\begin{document}
\maketitle

\begin{abstract}
We shall work with the so called duality triads following
Kwa\'sniewski. In particular in this note we propose some
extensions of them - hence we choose such  special class of triads
that admit - all at once - a unified combinatorial interpretation
in a way Konvalina does. The proposed extensions contain also
definition of operator valued arrays of Konvalina-like $C_{n,k}$ -
for the sake of future investigation and applications.
\end{abstract}


\section{Preliminaries.}
At first let us quote - following \cite{1,2} the basic definitions
introducing  the so called - duality triads.

As underlined in \cite{1,2} these duality triads arise as a
natural and sine qua non object for example  in dynamical data
basis` models. For this and other motivations and other "triad
sources" of inspiration - we refer to \cite{1,2} and references
therein. We shall consider now the infinite array of  numbers
  $c_{n,k}$ satisfying the following $3$-term recurrence equation
 (for $n,k$ being nonnegative integers and $i_k, q_k, d_k$ being integer numbers):
\begin{equation}\label{triadscoeff}
\left\{ {{\begin{array}{*{20}l}
{c_{n+1,k}=i_{k-1}c_{n,k-1} + q_{k}c_{n,k} +d_{k+1}c_{n,k+1}} \\
{c_{0,0}=1}\\
{c_{0,k}=0 \; \textrm{for}\; k>0}
\end{array}} } \right.
\end{equation}
Then the dual recurrence with respect to (\ref{triadscoeff}) is of
the form:
\begin{equation}\label{triadspoly}
\left\{ {{\begin{array}{*{20}l}
{x\phi_{k}(x)=d_k\phi_{k-1}(x)+q_k\phi_k(x) +i_k\phi_{k+1}(x)}\\
{\phi_0(x)=1}\\
{\phi_{-1}(x)=0}
\end{array}} } \right.
\end{equation}
where $\{\phi_k(x)\}_{k\geq 0}$ is a polynomial sequence such that
deg $\phi_k(x) =k$ and
\begin{equation}\label{triadsx}
x^{n}= \sum_{k=0}^{n}c_{n,k}\phi_k(x)
\end{equation}
Given the recurrence equations (\ref{triadscoeff}) and
(\ref{triadspoly}) one can derive the identity (\ref{triadsx})
(for proof see \cite{1}). From (\ref{triadsx}) one can observe
that the sequence of numbers $c_{n,k}$ is the sequence of
expansion coefficients of polynomials $x^n$ in the basis of
polynomials $\{\phi_k(x)\}_{k\geq 0}$. In combinatorics $c_{n,k}$
are called connection constants \cite{3,4}.

The equations
(\ref{triadscoeff}),(\ref{triadspoly}),(\ref{triadsx}) compose the
duality triads and the coefficients $c_{n,k}$ are called the triad
coefficients, while the polynomials $\{\phi_k(x)\}_{k\geq 0}$ -
are to be refer to as triad polynomials. The first important and
illustrative  examples of duality triad polynomials were delivered
in \cite{2} from where all our quotations come.

\section{Konvalina triads.}
In \cite{5,6} Konvalina presents  a unified simultaneous
combinatorial interpretation for both binomial coefficients of the
first and second kind and  Gaussian coefficients of the first and
second kind and Stirling numbers of the first and second kind.

Let us present it in brief. Consider as in \cite{6} $n$ different
boxes containing different balls - the $i$-th box contains $w_i$
balls ($w_i\geq 1$). The number $w_i$ is called the weight of box
$i$ and vector $\overrightarrow{w}=(w_1,w_2,\ldots,w_n)$ is the
weight vector.
\begin{defn}\label{K1}
The generalized Konvalina binomial coefficient of the first kind
with vector weight $\overrightarrow{w}$, denoted
$C_k^n(\overrightarrow{w})$, is the number of ways to select $k$
balls from $n$ different boxes $(i_1<i_2<\ldots<i_k)$ and then
taking one ball from each of the $k$ selected boxes (there are
$w_{i_1}w_{i_2}\cdots w_{i_k}$ possibilities):
\begin{equation}\label{Kon1}
C_k^n(\overrightarrow{w})=\sum_{1\leq i_1<i_2<\cdots <i_k\leq n}
w_{i_1}w_{i_2}\cdots w_{i_k}
\end{equation}
\end{defn}

\begin{defn}\label{K2}
The generalized Konvalina binomial coefficient of the second kind
with vector weight $\overrightarrow{w}$, denoted
$S_k^n(\overrightarrow{w})$, is the number of ways to select $k$
balls from $n$ not necessarily different boxes $(i_1\leq i_2\leq
\ldots \leq i_k)$ and then taking one ball from each of the $k$
selected boxes (there are $w_{i_1}w_{i_2}\cdots w_{i_k}$
possibilities):
\begin{equation}\label{Kon2}
S_k^n(\overrightarrow{w})=\sum_{1\leq i_1\leq i_2\leq \cdots \leq
i_k\leq n} w_{i_1}w_{i_2}\cdots w_{i_k}
\end{equation}
\end{defn}

Counting the ways of selections $k$ balls from $n$ different or
not necessarily different boxes, respectively, with dependence
either the last box was chosen or not one can derive the following
recurrences:
\begin{equation}\label{recKon1}
C_k^n(\overrightarrow{w})=C_k^{n-1}(\overrightarrow{w})+w_nC_{k-1}^{n-1}(\overrightarrow{w})
\end{equation}

\begin{equation}\label{recKon2}
S_k^n(\overrightarrow{w})=S_k^{n-1}(\overrightarrow{w})+w_nS_{k-1}^n(\overrightarrow{w})
\end{equation}
Setting
$\overline{S}_k^n(\overrightarrow{w})=S_{n-k}^k(\overrightarrow{w})$
one can get the equivalent to (\ref{recKon2}) recurrence:

\begin{equation}\label{recKon2'}
\overline{S}_k^n(\overrightarrow{w})= w_k
\overline{S}_k^{n-1}(\overrightarrow{w})+\overline{S}_{k-1}^{n-1}(\overrightarrow{w})
\end{equation}

\begin{obs}{\em Let $\overrightarrow{w}$ be a wieght vector. Then from Definitions (\ref{K1}),(\ref{K2})
one can observe:
\begin{enumerate}
\renewcommand{\labelenumi}{(\alph{enumi})}
\item $C_0^n(\overrightarrow{w})=1$ for $n\geq 0$ and $C_k^0(\overrightarrow{w})=0$ for $k>0$
\item $\overline{S}_0^n(\overrightarrow{w})=1$ for $n\geq 0$ and $\overline{S}_k^0(\overrightarrow{w})=0$ for $k>0$
\end{enumerate}}
\end{obs}

The Konvalina coefficients $C_k^n(\overrightarrow{w})$,
$\overline{S}_k^n(\overrightarrow{w})$ for the vector weight
$\overrightarrow{w}$ along with conditions given in the
observation above constitute two classes of triad coefficients
(\ref{triadscoeff}) with $i_{k-1}=w_n,\;q_k=1,\;d_k=0$ and
$i_k=1,\;q_k=w_k,\;d_k=0$, respectively. Then let us define two
classes of duality triads:

\begin{defn}\label{Kontriads1}
Konvalina triads of the first kind are triads satisfying the the
following two $2$-terms dual recurrence equations and the third
coming next equation interrelating the first two ones:
\begin{equation}\label{triadscoeffKon1}
\left\{ {{\begin{array}{*{20}l}
{C_{k}^{n+1}(\overrightarrow{w})=w_nC_{k-1}^n(\overrightarrow{w}) + C_{k}^n(\overrightarrow{w})} \\
{C_0^0(\overrightarrow{w})=1}\\
{C_k^0(\overrightarrow{w})=0 \; \textrm{for}\; k>0}
\end{array}} } \right.
\end{equation}

\begin{equation}\label{Kontriadspoly1}
\left\{ {{\begin{array}{*{20}l}
{x\phi_{k}(x)=\phi_k(x) +i_k\phi_{k+1}(x)\; \textrm{where}\; i_{k-1}=w_n}\\
{\phi_0(x)=1}\\
{\phi_{-1}(x)=0}
\end{array}} } \right.
\end{equation}
where $\{\phi_k(x)\}_{k\geq 0}$ is a polynomial sequence such that
deg $\phi_k(x) =k$ and
\begin{equation}\label{Kontriadsx1}
x^{n}= \sum_{k=0}^{n}C_{k}^n(\overrightarrow{w})\phi_k(x)
\end{equation}

\end{defn}

\begin{defn}\label{Kontriads2}
Konvalina triads of the second kind are triads satisfying the the
following two $2$-terms dual recurrence equations and the third
coming next equation interrelating the first two ones:
\begin{equation}\label{triadscoeffKon2}
\left\{ {{\begin{array}{*{20}l}
{\overline{S}_{k}^{n+1}(\overrightarrow{w})=\overline{S}_{k-1}^n(\overrightarrow{w}) +
w_k\overline{S}_{k}^n(\overrightarrow{w})} \\
{\overline{S}_0^0(\overrightarrow{w})=1}\\
{\overline{S}_k^0(\overrightarrow{w})=0 \; \textrm{for}\; k>0}
\end{array}} } \right.
\end{equation}

\begin{equation}\label{Kontriadspoly2}
\left\{ {{\begin{array}{*{20}l}
{x\phi_{k}(x)=w_k\phi_k(x) +\phi_{k+1}(x)}\\
{\phi_0(x)=1}\\
{\phi_{-1}(x)=0}
\end{array}} } \right.
\end{equation}
where $\{\phi_k(x)\}_{k\geq 0}$ is a polynomial sequence such that
deg $\phi_k(x) =k$ and
\begin{equation}\label{Kontriadsx2}
x^{n}=\sum_{k=0}^{n}\overline{S}_{k}^n(\overrightarrow{w})\phi_k(x)
\end{equation}
\end{defn}

Let us now quote from \cite{1,2} few simplest examples of triads
which are to be here now Konvalina triads of the first and second
kind.
\begin{ex} {\em Pascal triad:
$$
\left\{ {{\begin{array}{*{20}l}
{ \binom{n+1}{k}=\binom{n}{k-1} + \binom{n}{k} }\\
{ \binom{0}{0}=1 }\\
{\binom{k}{0}=0 \; \textrm{for}\; k>0}
\end{array}} } \right.
$$

$$
\left\{ {{\begin{array}{*{20}l}
{x\phi_{k}(x)=\phi_k(x) +\phi_{k+1}(x)}\\
{\phi_0(x)=1}\\
{\phi_{-1}(x)=0}
\end{array}} } \right.
$$ where
$$x^{n}=\sum_{k=0}^{n}\binom{n}{k}\phi_k(x).$$
Hence $$\phi_{k}(x)=(x-1)^k.$$ The coefficients
$c_{n,k}=\binom{n}{k}$ are Konvalina triads coefficients of the
first kind with wector weight $\overrightarrow{w}=(1,1,\ldots,1)$:
$$\binom{n}{k}=C_k^n(\overrightarrow{w}).$$}
\end{ex}

\begin{ex}{\em Stirling triad:
$$
\left\{ {{\begin{array}{*{20}l} { \left\{ {{\begin{array}{*{20}c}
{n}\\ {k}  \end{array}} } \right\} =
  k \left\{ {{\begin{array}{*{20}c} {n-1}\\ {k}  \end{array}} } \right\} +
  \left\{ {{\begin{array}{*{20}c} {n-1}\\ {k-1}  \end{array}} } \right\}} \\
{ \left\{ {{\begin{array}{*{20}c} {0}\\ {0}  \end{array}} } \right\}=1 }\\
{\left\{ {{\begin{array}{*{20}c} {0}\\ {k}  \end{array}} }
\right\}=0 \; \textrm{for}\; k>0}
\end{array}} } \right.
$$

$$
\left\{ {{\begin{array}{*{20}l}
{x\phi_{k}(x)=k\phi_k(x) +\phi_{k+1}(x)}\\
{\phi_0(x)=1}\\
{\phi_{-1}(x)=0}
\end{array}} } \right.
$$ where
$$x^{n}=\sum_{k=0}^{n}\left\{ {{\begin{array}{*{20}c} {n}\\ {k}  \end{array}} } \right\}\phi_k(x).$$
Hence $$\phi_{k}(x)=x^{\underline{k}}=x(x-1)(x-2)\cdots (x-k+1).$$
The coefficients $c_{n,k}=\left\{ {{\begin{array}{*{20}c} {n}\\
{k}  \end{array}} } \right\}$ are Konvalina triads coefficients of
the second kind with wector weight
$\overrightarrow{w}=(1,2,3,\ldots,n)$:
$$\left\{ {{\begin{array}{*{20}c} {n}\\ {k}  \end{array}} }
\right\}=S_{n-k}^k(\overrightarrow{w}).$$}
\end{ex}

\begin{ex} {\em $q$-Gaussian triad:
$$
\left\{ {{\begin{array}{*{20}l}
{ \binom{n+1}{k}_q=\binom{n}{k-1}_q + q^{k}\binom{n}{k}_q }\\
{ \binom{0}{0}_q=1 }\\
{\binom{k}{0}_q=0 \; \textrm{for}\; k>0}
\end{array}} } \right.
$$

$$
\left\{ {{\begin{array}{*{20}l}
{x\phi_{k}(x)=q^k \phi_k(x) +\phi_{k+1}(x)}\\
{\phi_0(x)=1}\\
{\phi_{-1}(x)=0}
\end{array}} } \right.
$$ where
$$x^{n}=\sum_{k=0}^{n}\binom{n}{k}_q\phi_k(x).$$
Hence $$\phi_{k}(x)=\prod_{i=0}^{k-1}(x-q^i).$$ The coefficients
$c_{n,k}=\binom{n}{k}_q$ are Konvalina triads coefficients of the
second kind with vector weight
$\overrightarrow{w}=(1,q,q^2,\ldots,q^{n-1})$:
$$\binom{n}{k}_q = S_k^{n-k+1}(\overrightarrow{w}).$$}
\end{ex}

Though the Stirling numbers of the first kind $\left[
{{\begin{array}{*{20}c} {n}\\ {k}  \end{array}} } \right]$
-satisfying  (\ref{recKon1}) -are the Konvalina generalized
binomial coefficients $\left(\left[ {{\begin{array}{*{20}c} {n}\\
{k}  \end{array}} } \right]=C_{n-k}^{n-1}(\overrightarrow{w})\;
\textrm{where}\; \overrightarrow{w}=(1,2,\cdots,n)\right)$ (see
\cite{5,6}) due  the  well known recurrence equation:
\begin{equation}\label{1}
\left\{ {{\begin{array}{*{20}l} { \left[ {{\begin{array}{*{20}c}
{n}\\ {k}  \end{array}} } \right] =
  (n-1) \left[ {{\begin{array}{*{20}c} {n-1}\\ {k}  \end{array}} } \right] +
  \left[ {{\begin{array}{*{20}c} {n-1}\\ {k-1}  \end{array}} } \right]} \\
{ \left[ {{\begin{array}{*{20}c} {0}\\ {0}  \end{array}} } \right]=1 }\\
{\left[ {{\begin{array}{*{20}c} {0}\\ {k}  \end{array}} }
\right]=0 \; \textrm{for}\; k>0}
\end{array}} } \right.
\end{equation}
they don't constitute the Konvalina triads coefficients
(\ref{triadscoeffKon1}). The non-triad polynomials that anyhow
might be eventually associated with the Stirling numbers of the
second kind - constitute - as seen from what follows - a one $n$ -
parameter family of polynomial sequences:

$$
\left\{ {{\begin{array}{*{20}l}
{x\phi_{k}(x)=(n-1) \phi_k(x) +\phi_{k+1}(x)}\\
{\phi_0(x)=1}\\
{\phi_{-1}(x)=0}
\end{array}} } \right. .
$$
The other important and illustrative  examples of duality triad
polynomials were delivered in \cite{1,2} from where all our
quotations come.

\section{Properties of Konvalina triads coefficients.}
The Konvalina triads coefficients deliver array of combinatoric`s
interpreted. These Konvalina triad coefficients  that share the
following elementary properties.

\begin{prop}\label{prop1}{\em
Let $\overrightarrow{w}$ be a vector weight of $n$ different
boxes. Then:
\begin{enumerate}
\item $C_{k}^{n}(\overrightarrow{w})=\sum_{i=0}^{n}w_iC_{k-1}^{i-1}(\overrightarrow{w})$
\item $S_{k}^{n}(\overrightarrow{w})=\sum_{i=0}^{n}w_iS_{k-1}^{i}(\overrightarrow{w})$
\end{enumerate}}
\end{prop}
\begin{proof} The proof runs counting ways of selections $k$ balls
from $n$ different boxes with vector weight $\overrightarrow{w}$
choosing at first $k$ different (not necessarily different) boxes,
respectively, and next taking one ball per box from selected boxes
with respect to the last box that has been selected. $A_i$ $(1\leq
i\leq n)$ denotes the class of $k$-selections that the box $i$ is
the last box has been selected. $|A_i|=w_iC_{k-1}^{i-1}$ as one
ball is taken from the box $i$ and $k-1$ balls from boxes
$\{1,2,3,\ldots,i-1\}$ (box repetitions not allowed).
$|A_i|=w_iS_{k-1}^{i}$ if box repetitions allowed (one ball is
taken from the box $i$ and $k-1$ balls from boxes
$\{1,2,3,\ldots,i\}$). Hence
$$C_{k}^{n}(\overrightarrow{w})=\sum_{i=1}^{n}w_iC_{k-1}^{i-1}(\overrightarrow{w}).$$
$$S_{k}^{n}(\overrightarrow{w})=\sum_{i=1}^{n}w_iS_{k-1}^{i}(\overrightarrow{w}).$$
\end{proof}

\begin{prop}\label{prop2}{\em
Let $\overrightarrow{w}$ be a vector weight of $n$ different
boxes. Then:
\begin{enumerate}
\item $C_{k}^{n+m}(\overrightarrow{w})=\sum_{i\geq
0}C_{i}^{n}(\overrightarrow{w})C_{k-i}^{m}(\overrightarrow{v})$
where $\overrightarrow{v}=(w_{n+1},w_{n+2},\ldots,w_{n+m})$.
\item $S_{k}^{n+m}(\overrightarrow{w})=\sum_{i\geq 0}S_{i}^{n}(\overrightarrow{w})S_{k-i}^{m}(\overrightarrow{v})$
where $\overrightarrow{v}=(w_{n+1},w_{n+2},\ldots,w_{n+m})$.
\end{enumerate}}
\end{prop}
\begin{proof}
Selecting $k$ balls from $n$ different boxes with vector weight
$\overrightarrow{w}$ choosing at first $k$ different (not
necessarily different) boxes, respectively, and next taking one
ball per box, $i$ balls are taken from $i$ boxes selected from
$\{1,2,\ldots,n\}$ and $k-i$ balls from $k-i$ boxes selected from
$\{n+1,n+2,\ldots,n+m\}$ $(0\leq i\leq k)$. Hence thesis.
\end{proof}

From Proposition (\ref{prop1}) one can get the following
well-known identities:
\begin{cor}{\em
Let us denote vectors:
$\overrightarrow{1}=(1,1,\ldots,1),\;\overrightarrow{i}=(1,2,3,\ldots,n),\;\overrightarrow{q}=(1,q,q^2,\ldots,q^{n-1})$.
\begin{enumerate}
\item $C_{k+1}^{n+k+1}(\overrightarrow{1})=\sum_{i=0}^{n+k+1}C_{k}^{i-1}(\overrightarrow{1})=\sum_{i=0}^{n}C_{k}^{k+i}(\overrightarrow{1})$
hence $$\binom{n+k+1}{k+1}=\sum_{i=0}^{n}\binom{k+i}{k}.$$
\item $C_{k+1}^{n+1}(\overrightarrow{1})=\sum_{i=0}^{n+1}C_{k}^{i-1}(\overrightarrow{1})=\sum_{i=0}^{n}C_{k}^{i}(\overrightarrow{1})$
hence $$\binom{n+1}{k+1}=\sum_{i=0}^{n}\binom{i}{k}.$$
\item $S_{k+1}^{n}(\overrightarrow{i})=\sum_{i=0}^{n}iS_{k}^{i}(\overrightarrow{i})=\sum_{i=0}^{n}iS_{k}^{i}(\overrightarrow{i})$
hence $$\left\{ {{\begin{array}{*{20}c} {n+k+1}\\ {n}
\end{array}}} \right\} = \sum_{i=0}^{n}i\left\{ {{\begin{array}{*{20}c} {k+i}\\ {i}  \end{array}} } \right\}.$$
\item $C_{k+1}^{n+k}(\overrightarrow{i})=\sum_{i=0}^{n+k}iC_{k}^{i-1}(\overrightarrow{i})=\sum_{i=k+1}^{n+k}iC_{k}^{i-1}(\overrightarrow{i})=\sum_{i=1}^{n}(k+i)C_{k}^{k+i-1}(\overrightarrow{i})$
hence $$\left( {{\begin{array}{*{20}c} {n+k+1}\\ {n}
\end{array}}} \right) = \sum_{i=0}^{n}(k+i)\left( {{\begin{array}{*{20}c} {k+i}\\ {i}  \end{array}} } \right).$$
\item $S_{k+1}^{n+1}(\overrightarrow{q})=\sum_{i=0}^{n+1}q^{i-1}S_{k}^{i}(\overrightarrow{q})=\sum_{i=0}^{n}q^iS_{k}^{i+1}(\overrightarrow{1})$
hence $$\binom{n+k+1}{k+1}_q=\sum_{i=0}^{n}\binom{k+i}{k}_q.$$
\end{enumerate}}
\end{cor}

Proposition (\ref{prop2}) in particular is the Cauchy identity:
\begin{multline*}
\binom{n+m}{k}=C_{k}^{n+m}(\overrightarrow{1})=\sum_{i\geq
0}C_i^n(\overrightarrow{1})C_{k-i}^m(\overrightarrow{1})=\sum_{i\geq
0}\binom{n}{i}\binom{m}{k-i}
\end{multline*}

\section{$\hat{q}_\psi$ Konvalina-like operators}
In this section we introduce as suggested by Kwasniewski to the
present author the operator valued  binomial array matrix elements
\cite{10,11} and then consequently operator valued Gaussian arrays
and operator valued Stirling arrays. Apart from the already
existing applications in \cite{10,11} a similar in character
operator valued infinite arrays appear also in physics (see for
more Katriel and Kibler \cite{Katriel}). The idea of considering
operator valued arrays (this time) in the setting of Konvalina
generalized binomial coefficients is presented in what follows.
All operators are supposed to act on the algebra of formal series
including the subalgebra $P=$\textbf{F}$[x]$ of polynomials of a
single variable $x$ over the field \textbf{F} of characteristic
$0$. We thus ensure the sufficient background for the
constructions to follow (see for more in \cite{12}).

\vspace{2mm}

 We shall use the upside-down notation introduced in \cite{7,8} and applied for example in
 \cite{10,11}
 which are of the source importance for what follows.
 At first then let us re-introduce this  $\psi$-notation.\\
 Consider $\Im $ - the family of functions` sequences such that:\\
$\Im = \{\psi;R \supset \left[ {a,b} \right]\;;\;q \in \left[
{a,b} \right]\;;\;\psi \left( {q} \right):Z \to F\;;\;\psi _{0}
\left( {q} \right) = 1\;;\;\psi _{n} \left( {q} \right) \ne
0;\;\psi _{ - n} \left( {q} \right)
= 0;\;n \in N\}$.\\

We shall call $\psi = \left\{ {\psi _{n} \left( {q} \right)}
\right\}_{n \ge 0} $ ; $\psi _{n} \left( {q} \right) \ne 0$; $n
\ge 0$ and $\psi _{0} \left( {q} \right) = 1$ an admissible
sequence. Consequently the symbol  $n_{\psi}  $ denotes
$$
n_{\psi}=\frac{\psi _{n - 1} \left( {q} \right)}{\psi _{n} \left(
{q} \right)},\;\;n \geq 0.$$

Then $\psi$-factorial and lower $\psi$-factorial are given
accordingly by
$$
n_{\psi}! \equiv \psi _{n}^{-1} \left( {q} \right) \equiv n_{\psi}
\left( {n - 1} \right)_{\psi} \left( {n - 2} \right)_{\psi} \left(
{n - 3} \right)_{\psi}  .... 2_{\psi} 1_{\psi}; \quad 0_{\psi}!=1
$$

 $n_{\psi} ^{\underline {k}}  = n_{\psi}  \left( {n - 1} \right)_{\psi}
\cdots \left( {n - k + 1} \right)_{\psi}$.

\begin{ex}\label{nq}{\em  Taking an admissible sequence $\psi_n(q)=(n_q!)^{-1}$
where $n_q=1+q+q^2+\cdots +q^{n-1}$ is well known $q$-deformation
of a natural number $n$ we obtain $n_{\psi}=n_q$. While $q=1$,
$n_q=n$.}
\end{ex}

For the definition of \textit{$q$-mutator} operator   to follow
see: Definition (3.2) in \cite{7} and (5.2) in \cite{8} for the
most general case.

\begin{defn} \cite{8,7} Let $\psi$ be an admissible sequence. The
$\hat{q}_\psi$ operator is a linear operator acting on the algebra
$P$ defined in the basis $\{x^n\}_{n\geq 0}$ as follows:
$$\hat{q}_\psi x^n=\frac{(n-1)_{\psi}-1}{n_{\psi}}x^n,\quad n\geq 0.$$
We shall call it after Kwa\'sniewski the $\hat{q}_\psi$-mutator
operator.
\end{defn}

\begin{defn} \cite{8,7} Let $\psi$ be an admissible sequence. The
$n_{\hat{q}_\psi}$ operator is a linear operator acting on the
algebra $P$ defined as follows:
$$n_{\hat{q}_\psi}=1+\hat{q}_\psi+\hat{q}_\psi^2+\cdots+\hat{q}_\psi^{n-1},\quad n>0$$
\end{defn}

\begin{obs}{\em
With the choice  $\psi_{n}(q)=(n_q)^{-1}$ one has of course:
\begin{enumerate}
\item $\hat{q}_\psi x^n=qx^n$ for $n\geq 0$. \item $\hat{q}_\psi^n
x^n=q^nx^n$ for $n\geq 0$. \item $n_{\hat{q}_\psi}x^n=n_qx^n$ and
for $q=1$ $n_{\hat{q}_\psi}x^n=nx^n$ is just the operator of
multiplication by $n$.
\end{enumerate}}
\end{obs}

Let us now come over to the announced extensions of
$\binom{n}{k}$, $\left[ {{\begin{array}{*{20}c} {n}\\ {k}
\end{array}} } \right]$ and $\left\{ {{\begin{array}{*{20}c} {n}\\
{k}  \end{array}} } \right\}$ as to become operator valued arrays
and  from now here on denoted correspondingly as :
$\binom{n}{k}_{\hat{q}_\psi}$, $\left[ {{\begin{array}{*{20}c}
{n}\\ {k}  \end{array}} }
\right]_{\hat{q}_\psi}$ and $\left\{ {{\begin{array}{*{20}c} {n}\\
{k}  \end{array}} } \right\}_{\hat{q}_\psi}$.

For that to do consider first the Konvalina binomial
coefficients.\\ Let
$\overrightarrow{w(\hat{q}_{\psi})}=(w_1(\hat{q}_{\psi}),
w_2(\hat{q}_{\psi}),\ldots,w_n(\hat{q}_{\psi}))$ denotes the
vector of linear operators such that
$w_i(\hat{q}_{\psi})=\sum_{k=0}^{i-1}a_{i,k}\hat{q}_{\psi}^k$
where $a_{i,k}\in$ {\textbf{F}} for $0\leq k\leq i-1$.

\begin{defn}
The $\hat{q}_{\psi}$ KonKwa   operator of the first kind, denoted
$C_{k}^{n}(\overrightarrow{w(\hat{q}_{\psi})})$ - is a linear
operator defined   on $P$  as follows:
$$C_k^n(\overrightarrow{w(\hat{q}_{\psi})})=\sum_{1\leq i_1<i_2<\cdots <i_k\leq n}
w_{i_1}(\hat{q}_{\psi})w_{i_2}(\hat{q}_{\psi})\cdots
w_{i_k}(\hat{q}_{\psi}.)$$
\end{defn}

\begin{defn}
The $\hat{q}_{\psi}$  KonKwa  operator of the second kind -
denoted $S_{k}^{n}(\overrightarrow{w(\hat{q}_{\psi})})$ - is a
linear operator  defined  on $P$  as follows:
$$S_k^n(\overrightarrow{w(\hat{q}_{\psi})})=\sum_{1\leq i_1\leq i_2\leq \cdots \leq i_k\leq n}
w_{i_1}(\hat{q}_{\psi})w_{i_2}(\hat{q}_{\psi})\cdots
w_{i_k}(\hat{q}_{\psi}).$$
\end{defn}

While setting
$\overrightarrow{w(\hat{q}_{\psi})}=(1,\hat{q}_{\psi},\hat{q}_{\psi}^2,\ldots,\hat{q}_{\psi}^{n-1})$
($1$ is identity operator) we obtain an extension of binomial
coefficient:
\begin{equation}
\binom{n}{k}_{\hat{q}_\psi}=S_k^{n-k+1}(\overrightarrow{w(\hat{q}_{\psi})})
\end{equation}
and choosing the $\psi$ admissible sequence
$\psi_n(q)=\frac{1}{n_q!}\quad
\binom{n}{k}_{\hat{q}_\psi}=\binom{n}{k}_{q}$ and putting  $q=1$
one ends up with  the ordinary binomial coefficient (treated as
linear operator of multiplication the polynomial $x^n$ by
$\binom{n}{k}_{q}$ and $\binom{n}{k}$, respectively).\\

Setting now
$\overrightarrow{w(\hat{q}_{\psi})}=(1_{\hat{q}_\psi},2_{\hat{q}_{\psi}},3_{\hat{q}_{\psi}},\ldots,n_{\hat{q}_{\psi}})$
we obtain an extension of Stirling numbers:

\begin{equation}
\left\{ {{\begin{array}{*{20}c} {n}\\ {k}
\end{array}}} \right\}_{\hat{q}_\psi}=S_{n-k}^{k}(\overrightarrow{w(\hat{q}_{\psi})})
\end{equation}

\begin{equation}
\left[ {{\begin{array}{*{20}c} {n}\\ {k}
\end{array}}} \right]_{\hat{q}_\psi}=C_{n-k}^{n-1}(\overrightarrow{w(\hat{q}_{\psi})})
\end{equation}

and again choosing the $\psi$ admissible sequence
$\psi_n(q)=\frac{1}{n_q!}\quad \left\{ {{\begin{array}{*{20}c}
{n}\\ {k}
\end{array}}} \right\}_{\hat{q}_\psi}=\left\{ {{\begin{array}{*{20}c} {n}\\ {k}
\end{array}}} \right\}_{q}$ and $\left[ {{\begin{array}{*{20}c} {n}\\ {k}
\end{array}}} \right]_{\hat{q}_\psi}=\left[ {{\begin{array}{*{20}c} {n}\\ {k}
\end{array}}} \right]_q$ and putting  $q=1$
one ends up with  the ordinary Stirling numbers (treated as linear
operator of multiplication  by $\left\{ {{\begin{array}{*{20}c}
{n}\\ {k}
\end{array}}} \right\}_{q}$, $\left\{ {{\begin{array}{*{20}c} {n}\\ {k}
\end{array}}} \right\}$, $\left[ {{\begin{array}{*{20}c} {n}\\ {k}
\end{array}}} \right]_{q}$, $\left[ {{\begin{array}{*{20}c} {n}\\ {k}
\end{array}}} \right]$ respectively).

The clue observation and hope for further investigation and
possible applications of the objects introduced above is that the
$\hat{q}_\psi$ KonKwa operators satisfy the same shape identities
and recurrence equations as the Konvalina generalized binomial
coefficients. These are all  obtained via replacing
$\overrightarrow{w}$ by $\overrightarrow{w(\hat{q}_\psi)}$.

\section{$\psi$-extensions of Konvalina binomial coefficients}
Another kind of extensions of Konvalina binomial coefficients one
can obtain by extending the vector weight of boxes. Let
$\overrightarrow{w_{\psi}}=\left(w_1,w_2,...,w_n \right)$ be such
extension of the vector weight $\overrightarrow{w}$ that $w_i \in
\{n_{\psi}:\;\;n \geq 0\}$. Then
\begin{equation}\label{psiKon1}
C_k^n(\overrightarrow{w_ \psi})=\sum_{1\leq i_1<i_2<\cdots
<i_k\leq n} w_{i_1}w_{i_2}\cdots w_{i_k}
\end{equation}
is called the $\psi$-extended generalized Konvalina binomial
coefficients of the first kind and
\begin{equation}\label{psiKon2}
S_k^n(\overrightarrow{w_\psi})=\sum_{1\leq i_1\leq i_2\leq \cdots
\leq i_k\leq n} w_{i_1}w_{i_2}\cdots w_{i_k}
\end{equation}
is called the $\psi$-extended generalized Konvalina binomial
coefficients of the second kind. The initial conditions are
following:\\
$S_0^n = C_0^n(\overrightarrow{w_\psi}) = 1$ for $n\geq 0$ and
$S_k^0 = C_k^0(\overrightarrow{w_\psi}) = 0$ for $k>0$.

Expectedly, these $\psi$-extended generalized Konvalina binomial
coefficients have combinatorial interpretation  for $n\geq 0$ and
$n_\psi$ being nonnegative integers - for example
Fibbonacci numbers.\\
Leave however this not easy question apart and let us do what we
can right now. Let us take a look at the vector weight
$\overrightarrow{w_\psi}= \left(1_\psi, 2_\psi, \ldots,n_\psi
\right)$. We obtain $\psi$-Stirling numbers:
\begin{equation} \label{psis}
\left[ {{\begin{array}{*{20}c} {n}\\ {k}  \end{array}} }
\right]_\psi = C_{n-k}^{n-1}(\overrightarrow{w_\psi}) =
\sum_{1\leq i_1< i_2< \cdots < i_{n-k}\leq n-1} (i_1)_\psi
(i_2)_\psi \cdots (i_{n-k})_\psi
\end{equation}
\begin{equation} \label{psiS}
\left\{ {{\begin{array}{*{20}c} {n}\\ {k}  \end{array}} }
\right\}_\psi = S_{n-k}^k(\overrightarrow{w_\psi})= \sum_{1\leq
i_1\leq i_2\leq \cdots \leq i_{n-k}\leq k} (i_1)_\psi (i_2)_\psi
\cdots (i_{n-k})_\psi
\end{equation}
Choosing the $\psi$ admissible sequence $\psi_n(q)=\frac{1}{n_q!}$
$\psi$-Stirling numbers become $q$-Stirling numbers $\left\{
{{\begin{array}{*{20}c} {n}\\{k}
\end{array}}} \right\}_{\psi}=\left\{ {{\begin{array}{*{20}c} {n}\\ {k}
\end{array}}} \right\}_{q}$, $\left[ {{\begin{array}{*{20}c} {n}\\ {k}
\end{array}}} \right]_{\psi}=\left[ {{\begin{array}{*{20}c} {n}\\ {k}
\end{array}}} \right]_q$ and for $q=1$ they are ordinary Stirling numbers.

One can easy observe that
\begin{equation} \label{genS}
\sum_{n\geq 0}\left\{ {{\begin{array}{*{20}c} {n}\\
{k} \end{array}} } \right\}_\psi x^n = \frac{x^k}{\left( 1-1_\psi
x\right)\left( 1-2_\psi x\right)\cdots \left( 1-k_\psi x\right)}
\end{equation}
and
\begin{equation} \label{gens}
\sum_{k\geq 0} \left[ {{\begin{array}{*{20}c} {n}\\ {k}
\end{array}} } \right]_\psi x^k = \Psi_{\overline{n}}(x),\quad
\Psi_{\overline{n}}(x)= x\left( x+1_\psi \right)\left( x+2_\psi
\right)\cdots \left( x+(n-1)_\psi \right)
\end{equation}
In \cite{13}  Kwa\'sniewski among others introduces following
Wagner \cite{14} the so called Comtet $\psi$-Stirling numbers of
the second kind, denoted there
$\left\{ {{\begin{array}{*{20}c} {n}\\
{k} \end{array}} } \right\}^\sim_\psi$ and defined equivalently
among others also by the identity (\ref{genS}) as well as the
Comtet $ \left[ {{\begin{array}{*{20}c} {n}\\
{k}\end{array}} } \right]^\sim_\psi$ Stirling numbers of the first
kind  different  from our  $ \left[ {{\begin{array}{*{20}c} {n}\\
{k}\end{array}} } \right]_\psi$ Stirling numbers of the first kind

Apart from  extension of Stirling numbers of the second kind to
Comtet numbers case there are also vastly  considered other
$q$-extended numbers
$S_q(n,k)=q^{\binom{k}{2}}\widetilde{S_q}(n,k)$ where
$\widetilde{S_q}(n,k)=\left\{ {{\begin{array}{*{20}c} {n}\\
{k} \end{array}} } \right\}^\sim_q$\quad (see [13] for
representative references and see also comments therein). May be
then for example these  $S_q(n,k)$ be also treated as another
$q$-extended generalized Konvalina binomial coefficients? And what
about their eventual operator valued arrays counterparts?\\

These and other questions  we leave for further examination and
investigation.\\

Extensions of generalized Konvalina binomial coefficients
introduced in this paper constitute only two classes of
extensions: $\hat{q}_\psi$-operator and $\psi$-extension Comtet one in terminology of \cite{13}.\\
The investigation is at preliminary stage - as seen from the all
above. We expect to say soon more on other extensions which can be
generated while choosing different vector weights
$\overrightarrow{w}$.

\end{document}